\documentclass[11pt,reqno]{amsart}
\usepackage{verbatim}
\usepackage{amssymb}
\usepackage{enumerate}
\usepackage[active]{srcltx}

\newif\ifdeveloping

\def\myheads#1;#2;{
\pagestyle{myheadings} \markboth{{\sc\hfill
#1\hfill\protect\makebox[0cm][r]{\rm\today}}}
{{\sc\protect\makebox[0cm][l]{\rm\today}\hfill #2\hfill}} }

\newcommand{\pcal}{{\mathcal P}}

\newcommand{\tcal}{{\mathcal T}}

\newcommand{\setm}{\setminus}
\newcommand{\empt}{\emptyset}
\newcommand{\subs}{\subset}

\def\<{\left\langle}
\def\>{\right\rangle}
\def\cf{\operatorname{cf}}

\def\OO{{\omega}}
\def\br#1;#2;{\bigl[ {#1} \bigr]^ {#2} }

\def\ooseq#1;#2;{\< {#1}_{#2}:{#2}<\oo\>}
\def\ooset#1;#2;{\{ {#1}_{#2}:{#2}<\oo\}}
\def\seq#1;#2;#3;{\< {#1}_{#2}:{#2}<#3\>}
\def\set#1;#2;#3;{\{ {#1}_{#2}:{#2}<#3\}}
\def\oseq#1;#2;{\< {#1}_{#2}:{#2}<\OO\>}
\def\oset#1;#2;{\{ {#1}_{#2}:{#2}<\OO\}}
\def\oosequ#1;#2;{\< {#1}^{#2}:{#2}<\oo\>}
\def\oosetu#1;#2;{\{ {#1}^{#2}:{#2}<\oo\}}
\def\sequ#1;#2;#3;{\< {#1}^{#2}:{#2}<#3\>}
\def\setu#1;#2;#3;{\{ {#1}^{#2}:{#2}<#3\}}
\def\osequ#1;#2;{\< {#1}^{#2}:{#2}<\OO\>}
\def\osetu#1;#2;{\{ {#1}^{#2}:{#2}<\OO\}}

\newcommand{\al}{\alpha}
\newcommand{\be}{\beta}

\def\to{\longrightarrow}

 \ifdeveloping
\usepackage[notref,notcite]{showkeys}
\fi
\usepackage{enumerate}

\newcommand{\prtime}{{\count0=\time\divide\count0 by 60
\count1=-\count0\multiply\count1 by 60 \advance\count1 by \time
\the\count0:\the\count1} }

\def\myheads#1;#2;{
\pagestyle{myheadings} \markboth{{\sc\hfill
#1\hfill\protect\makebox[0cm][r]{\rm\today; \prtime}}}
{{\sc\protect\makebox[0cm][l]{\rm\today;\ \prtime}\hfill
#2\hfill}} \thispagestyle{myheadings} }

\newcommand{\pib}{\operatorname{\pi}_B}
\newcommand{\piz}{\operatorname{\pi}}

\newcommand{\blocks}{\mathbb B}
\newcommand{\und}{X}

\newcommand{\ifu}{\operatorname{i}}

\newcommand{\concat}{\mathop{{}^{\frown}\makebox[-3pt]{}}}

\newcommand{\de}{\delta}
\newcommand{\ka}{\kappa}
\newcommand{\la}{\lambda}

\newtheorem{theorem}{Theorem}[section]
\newtheorem{proposition}[theorem]{Proposition}
\newtheorem{definition}[theorem]{Definition}

\newtheorem{lemma}[theorem]{Lemma}

\newcommand{\sub}{\subset}
\newcommand{\ga}{\gamma}
\newcommand{\om}{\omega}
\def\<{\left\langle}
\def\>{\right\rangle}

\newcommand{\un}{\bigcup}

\title{ On cardinal sequences of length $< \omega_3$}

\author{Juan Carlos Mart\'{\i}nez and Lajos Soukup}

\begin{document}

\maketitle

\footnotetext[1] { \hspace*{5mm} {\em Keywords and phrases}. locally compact scattered space, superatomic Boolean algebra, cardinal sequence.}

\begin{abstract}
We prove the  following consistency result for cardinal sequences of length $< \om_3$: if GCH holds and $\la \geq \om_2$ is a regular cardinal,
then in some cardinal-preserving generic extension  $2^{\om} = \la$ and for every ordinal $\eta < \om_3$ and every sequence $f = \langle \ka_{\al} : \al < \eta \rangle$
 of infinite cardinals with $\ka_{\al}\leq \la$ for $\al < \eta$ and $\ka_{\al} = \om$ if $\mbox{cf}(\al) = \om_2$, we have that
$f$ is the cardinal sequence of some LCS space.

Also, we prove that  for every specific uncountable cardinal $\lambda$ it is relatively consistent with ZFC that for every $\al,\be < \om_3$ with $\mbox{cf}(\al) < \om_2$ there is an LCS space $Z$ such that $\mbox{CS}(Z) = \langle \omega \rangle_{\alpha}\concat \langle \lambda \rangle_{\beta}$.
\end{abstract}

\section{Introduction}

\vspace{2mm} By an {\em LCS space} we mean a locally compact, Hausdorff and scattered space.  Recall that for an LCS space $X$ and an ordinal $\alpha$,
 the {\em ${\alpha}^{th}$- Cantor-Bendixson level} of $X$ is defined  by $I_{\al}(X) =$ the set of isolated
points of  $X\setminus \bigcup\{I_{\beta}(X): \beta < \alpha \}$.
  We define {\em the height of} $X$ as $ht(X) =$ the least ordinal $\de$ such that $I_{\de}(X) = \emptyset$, and we define its {\em reduced height} as $ht(X)^- =$ the least ordinal $\de$ such that $I_{\de}(X)$  is finite. Clearly, one has $\mbox{ht}^-(X) \leq \mbox{ht}(X) \leq \mbox{ht}^-(X) + 1$. We define the {\em width} of $X$ by $\mbox{wd}(X) = \mbox{sup}\{|I_{\al}(X)|: \al < \mbox{ht}(X)\}$. And we define the {\em cardinal sequence} of $X$  as $\mbox{CS}(X) = \langle |I_{\al}(X)| : \al < \mbox{ht}^-(X)\rangle.$ If $\kappa$ is an infinite cardinal and $\alpha$ is an ordinal, we denote by $\langle \kappa \rangle_{\alpha}$  the cardinal sequence $\langle \ka_{\be} : \be < \al \rangle$ where $\ka_{\be} = \ka$ for $\be < \al$.
If $f$ and $g$ are sequences of cardinals, we denote by $f\concat g$ the concatenation of $f$ with $g$.

Many authors have studied the possible sequences of infinite cardinals that can
arise as the cardinal sequence of an LCS space (or, equivalently, of a superatomic Boolen algebra). We refer the reader to the survey papers \cite{ba} and \cite{ma2} for a wide list of results on cardinal sequences of LCS spaces as well as examples and basic facts. It was proved by
Juh\'asz and Weiss that if $f = \langle \ka_{\al} : \al < \omega_1 \rangle$ is a sequence
of infinite cardinals, then $f$ is the cardinal sequence of an LCS space iff $\ka_{\be}\leq \ka^{\omega}_{\al}$ for every $\al < \be < \omega_1$ (see \cite[Theorem 5]{jw}). However, this result can not be extended to cardinal sequences of length $\omega_1 + 1$, since it was shown by Baumgarter in \cite{bs} that in the Mitchell Model there is no LCS space $X$ with $\mbox{CS}(X) = \langle \omega_1 \rangle_{\omega_1}\concat \langle \omega_2 \rangle$. Also,  a characterization
under GCH for cardinal sequences of length $< \omega_2$ was obtained  by Juh\'asz, Soukup and Weiss in \cite{jsw}. However, no characterization is known for cardinal sequences of length $\omega_2$. Nevertheless, it was shown by Baumgartner and Shelah in \cite{bs} that it is relatively consistent with ZFC that there is an LCS space of width $\omega$ and height $\omega_2$. This result was improved by Soukup in \cite{so}, where it was shown that if GCH holds and $\la \geq \om_2$ is a regular cardinal, then in some cardinal-preserving generic extension $2^{\om} = \la$ and every sequence $f = \langle \ka_{\al} : \al < \om_2 \rangle$ of infinite cardinals with $\ka_{\al} \leq \la$  is the cardinal sequence of some LCS space. However, the following basic proposition shows that Soukup's theorem can not be extended to cardinal sequences of length $<  \om_3$.

\begin{proposition}
Assume that $\kappa$ is a regular cardinal, $\eta > \kappa^{++}$ is an ordinal and $f = \langle \ka_{\xi} : \xi < \eta \rangle$ is the cardinal sequence of some LCS space. Assume that $\al < \eta$ is an ordinal with $\mbox{cf}(\al) > \kappa^+$ such that there is a strictly increasing sequence of ordinals $\langle \al_{\xi}: \xi < \mbox{cf}(\al) \rangle$ converging to $\al$ in such a way  that  $\ka_{\al_{\xi}} \leq \kappa$ for $\xi < \mbox{cf}(\al)$. Then, $\kappa_{\al}\leq \kappa$.
\end{proposition}

\vspace{2mm}
\noindent {\bf Proof}. Assume on the contrary that $\kappa_{\al} > \kappa$ . Let $X$ be an LCS space such that $\mbox{CS}(X) = f$. So, $|I_{\al}(X)| = \kappa_{\al}$. Let $Y$ be a subset of $I_{\al}(X)$ such that $|Y| = \kappa^+$. For every $x\in Y$ let $U_x$ be a compact open neighbourhood of $x$ such that $U_x \cap \un\{I_{\be}(X) : \al \leq \be < \eta \} = \{x\}$.  Clearly, for every $x,y\in Y$ with $x\neq y$, we have that $U_x\cap U_y \subset \un \{I_{\ga}(X) : \ga < \be \}$ for some $\be < \al$. Then, we define the function $F : [Y]^2 \longrightarrow \{\al_{\xi} : \xi < \mbox{cf}(\al) \}$ as follows. If $\{x,y\}\in [Y]^2$ we put

$$ F\{x,y\} = \mbox{ the least ordinal } \zeta < \mbox{cf}(\al) \mbox{ such that } U_x\, \cap \, U_y \subset \un \{I_{\ga}(X) : \ga < \al_{\zeta} \}.$$

\noindent Since $\mbox{cf}(\al) > \kappa^+$, there is a $\gamma <  \mbox{cf}(\al)$ such that for every $\{x,y\}\in [Y]^2$ we have $F\{x,y\} < \gamma$. But then $|I_{\al_{\ga}}(X)| = \kappa^+$, which contradicts the assumption that $\ka_{\al_{\xi}} \leq \kappa$ for $\xi < \mbox{cf}(\al)$. $\square$

\vspace{1mm} In particular, there is no LCS space $X$ with $\mbox{CS}(X) = \langle \omega \rangle_{\omega_2}\concat \langle \omega_1 \rangle$.

\vspace{1mm}
Then, we will show in Theorem 2.1 the following consistency result for cardinal sequences of length $< \om_3$: if GCH holds and $\la \geq \om_2$ is a regular cardinal,
then in some cardinal-preserving generic extension  $2^{\om} = \la$ and for every ordinal $\eta < \om_3$ and every sequence $f = \langle \ka_{\al} : \al < \eta \rangle$
 of infinite cardinals with $\ka_{\al}\leq \la$ for $\al < \eta$ and $\ka_{\al} = \om$ if $\mbox{cf}(\al) = \om_2$, we have that
$f$ is the cardinal sequence of some LCS space.

\vspace{1mm}
Also, we will prove in Theorem 3.1 that for every specific uncountable cardinal $\lambda$ it is relatively consistent with ZFC that for every $\al,\be < \om_3$ with $\mbox{cf}(\al) < \om_2$ there is an LCS space $Z$ such that $\mbox{CS}(Z) =
\langle \omega \rangle_{\alpha}\concat \langle
\lambda \rangle_{\beta}$.  This theorem improves the results shown in \cite{ro} and \cite[Section 2]{ms}.

\vspace{2mm} If $ T = \bigcup\{\{\al\}\times A_{\al} : \al < \eta \}$ where $\eta$ is a non-zero ordinal and each $A_{\al}$ is a non-empty set of ordinals,  then for every $s = \langle \al,\zeta\rangle\in T$ we write $\pi(s)= \al$ and $\rho(s) = \zeta$.

\vspace{2mm}
The following notion, which permits us to construct in a direct way LCS spaces  from  partial orders, will be used in our constructions.

\begin{definition}
 {\em  We say that $\tcal = \langle T,\preceq,i \rangle$ is an  {\em LCS poset}, if the following conditions hold:

\begin{enumerate}
\item $\langle T,\preceq \rangle$ is a partial order with $T = \bigcup \{T_{\al}: \al < \eta \}$ for some non-zero ordinal $\eta$
such that each $T_{\al} = \{\al\}\times A_{\al}$ where $A_{\al}$ is a non-empty set of ordinals.

\item  If $s \prec t$ then $\pi(s) < \pi(t)$.

\item If $\al < \be < \eta$ and $t\in T_{\be}$, then $\{s\in T_{\al} : s \prec t \}$ is infinite.

\item   $i : [T]^2 \rightarrow [T]^{<\om}$ such that for every $\{s,t\}\in [T]^2$ the following holds:

\begin{enumerate}[(a)]
\item If $v\in i\{s,t\}$, then $v\preceq s,t$.

\item If $u\preceq s,t$, then there is a $v\in i\{s,t\}$
such that $u\preceq v$.
\end{enumerate}

\end{enumerate}}

\end{definition}

\vspace{1mm}
If $\tcal = \langle T,\preceq,i \rangle$ is an LCS poset with $T = \bigcup \{T_{\al}: \al < \eta \}$ , we define its {\em associated LCS space} $X = X(\tcal)$ as follows. The underlying set of $X(\tcal)$ is $T$. If $x\in T$, we write  $C(x)= \{y\in T : y\preceq x\}$. Then, for every $x\in T$ we define a basic neighbourhood of $x$ in X as a set of the form $C(x)\setminus (C(x_1)\cup \dots \cup C(x_n))$ where $n < \omega$ and $x_1,\dots,x_n \prec x$. It can be checked that $X$ is a locally compact, Hausdorff, scattered space of height $\eta$ such that $I_{\al}(X) = T_{\al}$ for every $\al < \eta$ (see \cite{ba} for a proof). Then, we will say that $\langle |T_{\al}| : \al < \eta \rangle$ is the {\em cardinal sequence} of $\tcal$.

\vspace{1mm}
If $\tcal = \langle T,\preceq,i \rangle$ is an LCS poset and $S\subset T$ such that $i\{s,t\}\subset S$ for all $\{s,t\}\in [S]^2$, we define the {\em restriction of} $\tcal$ to $S$ as $\tcal\upharpoonright S = \langle S, \preceq \upharpoonright (S\times S), \\ i\upharpoonright [S]^2\rangle$.

\vspace{1mm}
The following notion, which is a refinement of the notion of a skeleton  given in  \cite[Definition 1.5]{rv}, will also  be
needed to show our results.

\begin{definition}{\em Assume that $\tcal = \<T,\preceq,i\>$ is an LCS poset
with $T =  \bigcup \{\{\al\}\times A_{\al}: \al < \eta \}$. Let $f$ be the cardinal sequence of $\tcal$.
Then, we say that $\tcal$ is an $f$-{\em skeleton}, if for every $\al <  \eta$ there is a countable subset
$O_{\al}\in [A_{\al}]^{\om}$ such that $s \prec t$ and $\pi(s) = \al$ implies $\rho(s)\in O_{\al}$, and in such a way that  the following two conditions hold:

\begin{enumerate}

\item If $\al < \eta$ and $s,t\in \{\al\}\times O_{\al}$  with  $\rho(s)\neq \rho(t)$, then $i\{s,t\} = \empt$.

\item If $\al + 1 < \eta$, $t\in \{\al + 1\}\times A_{\al + 1}$ and $s \prec t$, then there is a $u\in \{\al\}\times O_{\al}$ such that $s\preceq u \prec t$.
\end{enumerate}}

\end{definition}

\section{ A general consistency result}

\vspace{2mm}
 In this section, our aim is to prove the following result.

\begin{theorem}
If GCH holds and $\la \geq \om_2$ is a regular cardinal,
then in some cardinal-preserving generic extension  $2^{\om} = \la$ and for every ordinal $\eta < \om_3$ and every sequence $f = \langle \ka_{\al} : \al < \eta \rangle$
 of infinite cardinals with $\ka_{\al}\leq \la$ for $\al < \eta$ and $\ka_{\al} = \om$ if $\mbox{cf}(\al) = \om_2$, we have that
$f$ is the cardinal sequence of some LCS space.
\end{theorem}

In order to prove Theorem 2.1, we will use the following refinement of the notion of Shelah's $\Delta$-function due to Soukup.

\begin{definition}
{\em
 A function $f : [\om_2 \times \lambda]^2 \rightarrow [\om_2]^{< \om}$ is a  $\Delta$($\om_2 \times \la)$-{\em function}, if $f\{x,y\}\sub \mbox{min}\{\pi(x),\pi(y)\}$ for each $\{x,y\}\in [\om_2 \times \lambda]^2$ and for every uncountable subset $\{d_{\al} : \al < \om_1 \}$ of $[\om_2 \times \la]^{< \om}$ there are ordinals $\al < \be < \om_1$ such that for every $x\in d_{\al}\setm d_{\be}$, $y\in d_{\be} \setm d_{\al}$ and $z\in d_{\al}\cap d_{\be}\cap (\om_2\times\om)$ the following conditions hold:

\begin{enumerate}
\item if $\pi(z) < \pi(x),\pi(y)$ then $\pi(z)\in f\{x,y\}$,

\item  if $\pi(z) < \pi(y)$ then $f\{x,z\}\sub f\{x,y\}$,

\item  if $\pi(z) < \pi(x)$  then $f\{y,z\}\sub f\{x,y\}$.
\end{enumerate}}

\end{definition}

The following result was shown in \cite{so}.

\begin{theorem}
If GCH holds and $\la \geq \om_2$ is a regular cardinal, then in some cardinal-preserving generic extension $\lambda$ is a regular cardinal with $\la^{\om_1} = \la$ and
there is a $\Delta$($\om_2 \times \la)$-function.
\end{theorem}

\vspace{1mm}
 The following result can be proved by means of a slight refinement of the arguments given in \cite[Theorem 1.12]{rv} and \cite[Theorem 4.2]{so}.

\begin{theorem}
Assume that $\la \geq \om_2$ is a regular cardinal such that $\la^{\om} = \la$. Assume that there is a $\Delta(\om_2\times \la)$-function. Then,
in some c.c.c. generic extension $2^{\om} = \la$ and there is a $\langle \la \rangle_{\om_2}$--skeleton.
\end{theorem}

\begin{theorem}
Assume that $\la \geq \om_2$ is a regular cardinal with $2^{\om} = \la$, $\eta < \om_3$ is an ordinal  and there is
a $\langle \la \rangle_{\om_2}$-skeleton. Then, if $f = \langle \ka_{\al} : \al < \eta \rangle$ is a sequence of infinite cardinals with $\ka_{\al}\leq \la$
for $\al < \eta$ and $\ka_{\al} = \om$ if $\mbox{cf}(\al) = \om_2$, we have that $f$ is the cardinal sequence of some LCS space.
\end{theorem}

\vspace{2mm}
Note that Theorem 2.1 follows immediately from Theorems 2.3, 2.4 and 2.5.

\vspace{3mm}
\noindent {\bf Proof of Theorem 2.5}. Assume that $\la\geq\om_2$ is a regular cardinal such that $2^{\om} = \la$ and that there is a  $\langle \la \rangle_{\om_2}$-skeleton
$\tcal = \<T,\preceq,i\>$. We may assume that $O_{\al} = \om$ for every $\al < \om_2$. Proceeding by transfinite induction on $\eta < \om_3$, we show that if $f = \langle \ka_{\al}: \al < \eta \rangle$ is a sequence of
infinite cardinals with $\ka_{\al}\leq \la$
for $\al < \eta$ and $\ka_{\al} = \om$ if $\mbox{cf}(\al) = \om_2$, then there is an LCS space whose cardinal sequence is $f$.
If $\eta \leq \om_2$, we are done by the existence
of the  $\langle \la \rangle_{\om_2}$-skeleton  $\tcal$. So, assume that $\om_2 < \eta < \om_3$ and $f$ is as above.

 \vspace{1mm} First, suppose that $\eta$ is a limit ordinal. We assume that $\cf(\eta) = \omega_2$. Otherwise, the argument is similar. We distinguish two cases. First, suppose
 that there is a strictly increasing sequence of ordinals $\langle \gamma_{\xi}: \xi < \om_2\rangle$ converging to $\eta$ such that
$\mbox{cf}(\gamma_{\xi})= \om_2$ for every $\xi < \om_2$.  Note that if $\xi < \om_2$ is a limit ordinal, then
$\gamma_{\xi} > \mbox{sup} \{\gamma_{\mu} : \mu < \xi \}$. Now, we define the sequence of ordinals $\<\al_{\xi} : \xi < \om_2 \>$ as follows. We put $\al_0 = 0$. If $\xi$ is a successor ordinal,  we put $\al_{\xi} = \ga_{\xi}$. And if $\xi$ is a limit ordinal, we define $\al_{\xi} = \mbox{sup}\{\al_{\mu} : \mu < \xi \}$.

Now, for every $\xi < \om_2$, we
define $S_{\xi}$ as follows. We put $S_0 = \{0\}\times \ka_0$. If $\xi$ is a successor ordinal, we have $\ka_{\al_{\xi}} = \om$ and then we put  $S_{\xi} = \{\xi\}\times\om$. And if $\xi$ is a limit ordinal, we define $S_{\xi} = \{\xi\}\times \ka_{\al_{\xi}}$. We put $S = \bigcup\{S_{\xi} : \xi < \om_2 \}$. Note that $i\{s,t\}\subset \om_2 \times \om \subset  S$ for all $\{s,t\}\in [S]^2$. Then, we define ${\mathcal S} = \tcal \upharpoonright S = \<S,\preceq \upharpoonright (S\times S), i\upharpoonright [S]^2\>$. Clearly, ${\mathcal S}$ is an $h$-skeleton for $h = \langle \kappa_{\al_{\xi}} : \xi < \om_2 \rangle$.

 In order to carry out the construction, for every ordinal $\xi < \om_2$ we will insert an
adequate LCS space between the $\xi$-th and the $(\xi + 1)$-th level of the space associated with ${\mathcal S}$.
For every successor ordinal $\xi = \mu + 1 <  \om_2$, we define  $\delta_{\xi} = \mbox{o.t.}(\al_{\xi}\setm \al_{\mu})$.  By the induction hypothesis, for every successor ordinal $\xi = \mu + 1 < \om_2$ there is an LCS space whose cardinal sequence is $f_{\xi} = \langle \omega \rangle\concat \langle \ka_{\al_{\mu} + \zeta} : 0 < \zeta < \de_{\xi}\rangle$.

Let $S' = \bigcup \{ \{\xi\} \times \omega : \xi < \om_2 \mbox{ is a successor ordinal} \}$. For every successor ordinal $\xi = \mu + 1 < \om_2$ and every element $x\in \{\xi\} \times \omega$, we put $C_x = \{ y \in \{\mu\}\times \om : y \prec x  \}$. Note that since ${\mathcal T}$ is a skeleton, $C_x \cap C_y = \empt$ if $\pi(x) = \pi(y)$ and $x\neq y$.
Then, for every successor ordinal $\xi < \om_2$ and every element $x\in \{\xi\} \times \omega$, we consider a compact Hausdorff scattered space $Z_x$ of height $\delta_{\xi} + 1$ such that $\mbox{CS}(Z_x) = f_{\xi}$ and in such a way that $I_0(Z_x) = C_x$ and $I_{\delta_{\xi}}(Z_x) = \{x\}$. Also, we assume that for every $x\in S'$ we have $Z_x\cap S = C_x \cup \{x\}$ and that for every $x,y\in S'$ with $x\neq y$ we have $(Z_x \setminus \{x\}) \cap (Z_y \setminus \{y\}) = \emptyset$.

 Now, we  define the required space $Z$ as follows. Its underlying set is $S \cup \bigcup \{Z_x : x \in S'\}$.
For every $z\in S$ we put

$$W_z = \{y\in S : y \preceq z \}\cup \bigcup\{Z_y :  y \preceq z \mbox{ and } y\in S' \}.$$

\noindent Clearly, $x \prec z$ implies $W_x \subset W_z$. Also, if $z\in Z_x$ for some $x\in S'$ and $U$ is an
 open neighbourhood of $z$ in $Z_x$, we define

$$U^* = U \cup \bigcup\{W_y :   y\in C_x \cap U \}.$$

\noindent Note that as ${\mathcal T}$ is a skeleton, $W_x = Z^*_x$ for every $x\in S'$.

\vspace{1mm} Assume that $v\in S \cup \bigcup \{Z_x : x \in S'\}$.
Then, if $v\in S \setminus S'$  we define a basic neighbourhood of $v$ in $Z$ as a set of the form $W_v \setminus (W_{v_1}\cup\dots\cup  W_{v_n})$ where $n < \omega$ and $v_1,\dots,v_n \prec v$. If $v\in S'$, a basic neighbourhood of $v$ in $Z$ is a set $U^*$  where $U$ is a compact open neighbourhood  of $v$ in $Z_v$. Otherwise, we have that $v \in Z_x$ for a unique $x\in S'$, and then we  define a basic neighbourhood of $v$ in $Z$ as a set $U^*_v$ where $U_v$ is a compact open neighbourhood of $v$ in $Z_x$. In order to show that this collection of sets is in fact a base, suppose that $x\in W_y \setminus (W_{y_1}\cup \dots \cup W_{y_n})$ where $y_1,\dots , y_n \prec y \in S\setminus S'
$ and $x\in U^*_z$ where $z\in Z_w$ for some $w\in S'$ and $U_z$ is a compact open neighbourhood of $z$ in $Z_w$. The other cases are easier to verify. Without loss of generality, we may assume that $x\in U_z\setminus S$.  Since $x\in  (W_y \setminus (W_{y_1}\cup \dots \cup W_{y_n}))\cap  U_z$, we have that $w \preceq y$ and $w \not\preceq y_i$ for $i\in \{1,\dots , n\}$. Let $\{v_1,\dots,v_k\} = i\{w,y_1\}\cup \dots \cup i\{w,y_n\}$. Since ${\mathcal T}$ is a skeleton, there are $u_1,\dots , u_k \in C_w$ such that $v_i \preceq u_i$ for $i\in \{1,\dots, k\}$. Let $U_x$ be a compact open neighbourhood of $x$ in $Z_w$ such that $U_x\subset U_z$ and $U_x \cap \{u_1,\dots , u_k \} =  \emptyset$. Then, $U^*_x$ is as required. Also, it is easy to check that each element of this base is a clopen set.



We show that $Z$ is Hausdorff. Assume that $\{s,t\} \in [Z]^2$. We may assume that $s,t\not\in S$. Otherwise, the argument is easier.
Let $x,y\in S'$ be such that $s\in Z_x$ and  $t\in Z_y$. If $x = y$, the case is obvious. So, assume that $x\neq y$. Note that if $\pi(x) = \pi(y)$, $U_s$ is an open neighbourhood of $s$ in $Z_x$ and $U_t$ is an open neighbourhood of $t$ in $Z_y$, then as ${\mathcal T}$ is a skeleton we have  $U^*_s \cap U^*_t = \empt$.
So, suppose that $\pi(x) < \pi(y)$. If $x \prec y$, then $W_x\cap (W_y \setminus W_x) = \empt$ and clearly $W_x$ is a neighbourhood of $s$ and $W_y \setminus W_x$ is a neighbourhood of $t$.
 Now, assume that $x,y$ are $\preceq$-incomparable.  Put $i\{x,y\} = \{v_1,\dots,v_n\}$.  We show that $W_x \cap W_y \subset W_{v_1} \cup \dots \cup W_{v_n}$. For this, suppose that $u\in W_x \cap W_y$. Without loss of generality, we may assume that $u\in Z\setminus S$. Let $w\in S'$ be such that $u\in Z_w$. Since $u\in W_x\cap W_y$,  it follows that $w\prec x,y$, hence  $w\preceq v_i$ for some $i\in \{1,\dots,n\}$, and so $u\in W_{v_i}$. Therefore, $W_x \cap W_y \subset W_{v_1} \cup \dots \cup W_{v_n}$, and so we are done because $s\in W_x$ and $t\in W_y\setminus (W_{v_1} \cup \dots \cup W_{v_n})$.

Also, proceeding by transfinite induction on $\pi(x)$, we can  verify that $W_x$ is a compact neighbourhood of $x$ in $Z$ for every $x\in S$, and that if $z \in Z_x\setminus S$ for some $x\in S'$ and $U_z$ is a compact open neighbourhood of $z$ in $Z_x$ then $U^*_z$ is a compact neighbourhood of $z$ in $Z$. Therefore, $Z$ is locally compact.

 On the other hand, it is easy to see that $Z$ is a scattered space with $\mbox{CS}(Z) = f$.

\vspace{1mm}
Next, assume that there is no strictly increasing sequence of $\om_2$ many ordinals of
cofinality $\om_2$ converging to $\eta$. Let $\ga = \mbox{sup}\{\xi + 1 : \xi < \eta,
\mbox{cf}(\xi) = \om_2 \}$. Put $g = \langle \la_{\xi} : \xi < \eta \rangle$
where $\la_{\xi} = \om$ for $\xi\leq \ga$ and $\la_{\xi} = \ka_{\xi}$ for $\ga < \xi < \eta$. First, we construct an LCS space $Y$ with $\mbox{CS}(Y) = g$.
 Put $\zeta = \mbox{o.t.}(\eta\setminus \ga)$. Since there
is a $\langle \la \rangle_{\om_2}$--skeleton and in $\eta\setm \ga$ there is no ordinal of cofinality $\om_2$, it follows that there is an LCS space $X$ such that $\mbox{CS}(X) = \< \lambda_{\xi} : \xi < \zeta \>$ where $\lambda_0 = \om$ and $\la_{\xi} = \ka_{\ga + \xi}$
for $0 < \xi < \zeta$. Let $\{x_n : n\in \omega \}$ be an enumeration without repetitions of the elements of $I_0(X)$. By the induction hypothesis, for every $n < \om$ we can consider a compact Hausdorff scattered space $X_n$ of height $\gamma + 1$ such that $I_{\gamma}(X_n) = \{x_n\}$ and $\mbox{CS}(X_n) = \< \om \>_{\gamma}$ and in such a way that $X_n \cap X = \{x_n\}$ and $X_n \cap X_m = \empt$ for $n\neq m$. Then, the underlying set of $Y$ is $X \cup \bigcup \{X_n : n < \omega \}$. If $x\in X_n$ for some $n < \om$, a basic neighbourhood of $x$ in $Y$ is an open neighbourhood of $x$ in $X_n$. And if $x\in X\setminus I_0(X)$, then a basic neighbourhood of $x$ in $Y$ is a set of the form $U \cup \bigcup\{X_n : x_n\in U, n < \om \}$ where $U$ is an open neighbourhood of $x$ in $X$. Clearly, $Y$ is an LCS space with $\mbox{CS}(Y) = g$. Also, by the induction hypothesis, there is an LCS space $Z$ such that $\mbox{CS}(Z) = \< \kappa_{\xi} : \xi \leq \gamma \>$. We may assume that $Y\cap Z = \empt$. Then, the topological sum of $Y$ and $Z$ is the required LCS space of cardinal sequence $f$.

Finally, assume that $\eta = \ga + 1$ is a successor ordinal. First, suppose that $\ga$ is a limit
ordinal. If $\mbox{cf}(\ga) = \om$, then as $\ka_{\ga}\leq \la = 2^{\om}$, we can  carry out
the construction of the desired LCS space  by using an almost disjoint family of infinite subsets of
$\om$ of size $\ka_{\ga}$. If $\mbox{cf}(\ga) = \om_1$, we can proceed by means of an argument similar
to the one given in the case in which $\eta$ is a limit ordinal of cofinality $\om_2$.
And if $\mbox{cf}(\ga) = \om_2$, then $\ka_{\ga} = \om$ and so the construction is straightforward.
Now, suppose that $\ga = \de + 1$ is a successor ordinal. Since $\kappa_{\gamma} \leq 2^{\om}$, there is an LCS space $X$ such that $|I_0(X)| = \om$, $|I_1(X)| = \kappa_{\gamma}$ and $I_2(X) = \empt$. Then, proceeding as above, we can construct an LCS space $Y$ such that $X$ is a closed subspace of $Y$ and $\mbox{CS}(Y) = \< \lambda_{\xi} : \xi < \eta \>$ where $\lambda_{\xi} = \om$ for $\xi \leq \delta$ and $\lambda_{\gamma} = \kappa_{\gamma}$. Also, by the induction hypothesis, there is an LCS space $Z$ such that $\mbox{CS}(Z) = \< \kappa_{\xi} : \xi \leq \delta \>$ and $Y\cap Z = \empt$. Clearly, the topological sum of $Y$ and $Z$ is the required LCS space. $\square$

\section{A further construction of LCS spaces with countable levels}
\setcounter{case}{0}

In this section, our aim is to prove the following result.

\begin{theorem}
If $V = L$ holds and $\lambda$ is an uncountable cardinal, then in some cardinal-preserving generic extension we have that for every $\al,\be < \om_3$ with $\mbox{cf}(\al) < \om_2$ there is an LCS space $Z$ such that $\mbox{CS}(Z) =
\langle \omega \rangle_{\alpha}\concat \langle
\lambda \rangle_{\beta}$.
\end{theorem}

 In order to prove Theorem 3.1, we shall use the main result of \cite{ko}. First, we need to introduce the following notion of special function due to Koszmider.

 \begin{definition}
{\em  Assume that $\ka$ and $\la$ are infinite cardinals such that $\ka$ is
regular and $\ka <\la$. We say that a function $F :
[\la]^2\longrightarrow \ka^+$ is a $\ka^+$-{\em strongly unbounded
function on} $\la$, if  for every ordinal $\delta < \ka^+$, every
cardinal $\nu < \ka$ and every family $A\subset [\la]^{\nu}$ of
pairwise disjoint sets with $|A|=\ka^+$,
 there are different $a,b\in A$
such that $F\{\al,\be\}> \delta$ for every $\al\in a$ and $ \be\in b$.}  \end{definition}

The following result was proved in \cite{ko}.

\begin{theorem}
If $\ka$ and $\la$ are infinite cardinals such that $\ka^{+++} \leq
\la$, $\ka^{<\ka}=\ka$ and $2^{\ka}= \ka^+$, then in some cardinal-preserving generic extension
there is a $\ka^+$-strongly unbounded function on $\la$.
\end{theorem}

 In order to prove Theorem 3.1,  we need some preparation. Assume  that $\lambda$ is an uncountable cardinal. If $\gamma$ is an ordinal, we put

$$Y_{\gamma} = \gamma \times \omega \cup (\{\gamma, \gamma + 1\}\times \lambda).$$

 Let

$$ {\Bbb B } = \{S\} \cup  \lambda.$$

Let

$$ {\Bbb B }^{(\gamma)}_S = \gamma \times \omega$$

and

$$ {\Bbb B }^{(\gamma)}_{\zeta} = \{\gamma\}\times [\omega\cdot \zeta, \omega\cdot \zeta + \om) \cup \{\langle \gamma + 1,\zeta \rangle \}$$

\noindent for $\zeta \in \lambda$.

\vspace{4mm}
Clearly, $\{{\Bbb B }^{(\gamma)}_t : t\in {\Bbb B } \}$ is a partition of $Y_{\gamma}$. We define
\begin{equation*}
\pib:\und\to \blocks  \text{ by the formula } x\in \blocks^{(\gamma)}_{\pib
(x)}.
\end{equation*}

\vspace{1mm} Then, the following notion will be used in the proof of Theorem 3.1.

\begin{definition}

 {\em If $\gamma$ is an ordinal, we say that   $\tcal = \<T,\preceq,i\>$ is an {\em adequate} $\gamma$-{\em poset}, if $\tcal$ is an LCS poset with $T = Y_{\gamma}$ such that for every  $\zeta\in \lambda$ the following conditions hold:

 \begin{enumerate}
 \item If $x\in \blocks^{(\gamma)}_{\zeta}$ with $\pi(x) = \gamma$, then $x \prec \langle \gamma + 1, \zeta \rangle$ and $x \not\prec \langle \gamma + 1, \xi \rangle$ for $\xi\neq \zeta$.
 \item The restriction of $\tcal$ to ${\Bbb B}^{(\gamma)}_S \cup {\Bbb B}^{(\gamma)}_{\zeta}$ is an $\langle \omega \rangle_{\gamma + 1}\concat \langle 1 \rangle$-skeleton.

 \end{enumerate}}

  \end{definition}

  \begin{proposition}
  If there is an adequate $\omega_1$-poset, then there are an  adequate $\omega$-poset and an  adequate $1$-poset.
  \end{proposition}

  \noindent {\bf Proof}. Assume that there is an adequate $\omega_1$-poset. Since its cardinal sequence is $\langle \omega \rangle_{\omega_1}\concat \langle
\lambda \rangle_{2}$, we have that $2^{\omega}\geq \lambda$. Then, we construct an adequate $\omega$-poset $\tcal = \langle Y_{\om},\preceq,i \rangle$ as follows. First, for every $n < \om$, we consider an $\langle \omega \rangle_{n}\concat \langle 1\rangle$-skeleton  $\tcal_n = \langle T_n,\preceq_n,i_n \rangle$ such that $\{T_n : n < \omega \}$ is a partition of $\om\times \om$. Also, we assume that for each $n < \om$, there is a top point $v_n$ in $T_n$ such that $u <_n v_n$ for every $u\in T_n\setminus \{v_n\}$. Now, let $Y = \{\om\}\times \lambda$ and  $Y' = \{\om + 1\}\times \lambda$. For $\xi < \lambda$, put $y'_{\xi} = \langle\om + 1, \xi\rangle$. Since  $2^{\omega}\geq \lambda$, there is a pairwise disjoint family $\{a_{\xi} : \xi < \lambda\}$ of infinite subsets of $\omega$. Then, for every $\{x,y\}\in [(\om \times \om)\cup Y']^2$, we put $x \prec' y$ iff either $x\prec_n y$ for some $n < \omega$ or $x\in T_n$, $n\in a_{\zeta}$ and $y = y'_{\zeta}$.

 Now, for $\nu < \lambda$ and $n < \omega$, we put $z_{\nu,n} = \langle \om, \om\cdot\nu + n \rangle$. For $\xi < \lambda$, let $\{a_{\xi,n} : n < \om \}$ be a partition of $a_{\xi}$ into infinite subsets. Then, for every $\{x,y\}\in [Y_{\om}]^2$, we put $x \prec y$ iff one of the following conditions holds:

 \begin{enumerate}[(a)]
 \item  $x\prec' y$,

 \item for some $\nu < \lambda$, $x\in T_n$ for some $n\in a_{\nu,n}$ and $y = z_{\nu,n}$,

 \item for some $\nu < \lambda$ and $n < \omega$, $x = z_{\nu,n}$ and $y = y'_{\nu}$
 \end{enumerate}

 Now,  we define the infimum function $i$ as follows. Assume that $x,y\in Y_{\om}$ are $\preceq$-incomparable. If $x,y\in T_n$ for some $n\in \omega$, we put $i\{x,y\} = i_n\{x,y\}$. If $x,y\in Y \cup Y'$, we define $i\{x,y\} = \{v : v = v_n$ for some $ n < \om, v\prec x$ and $v \prec y\}$, which is a finite set because $\{a_{\xi} : \xi < \lambda\}$ is pairwise disjoint.
 And we put $i\{x,y\} = \emptyset$ otherwise. It is easy to check that $\langle Y_{\om},\preceq,i \rangle$ is as required. And by means of a simpler argument, one can show that there is an adequate $1$-poset.
$\square$

\vspace{1mm} In order to prove Theorem 3.1, we need to show the following lemma.

\begin{lemma} Assume that $\square_{\omega_1}$ holds and there is an adequate $\om_1$-poset. Then, in some cardinal-preserving generic extension we have that for every $\al,\be < \om_3$ with $\mbox{cf}(\al) < \om_2$ there is an LCS space $Z$ such that $\mbox{CS}(Z) =
\langle \omega \rangle_{\alpha}\concat \langle
\lambda \rangle_{\beta}$.

\end{lemma}

\noindent {\bf Proof}.  Since $\square_{\omega_1}$ holds, there is a cardinal-preserving generic extension $N$ where there is an LCS poset of width $\om$ and height $\eta$ for every ordinal $\eta < \om_3$ (see \cite{rv} or \cite{ma1}). By using Proposition 3.5, we have that in $N$ there is an adequate $\gamma$-poset for $\gamma \in \{1,\omega,\omega_1 \}$. Then, we prove that in $N$ there is an LCS space $Z$ such that $\mbox{CS}(Z) = \langle \omega \rangle_{\alpha}\concat \langle \lambda \rangle_{\beta}$ for every $\al,\be < \om_3$ with $\mbox{cf}(\al) < \om_2$ . Without loss of generality, we may assume that $\be$ is an infinite successor ordinal $\be' + 1$.  Assume that $\mbox{cf}(\al) = \om_1$. Let  $\tcal = \<T,\preceq,i\>$  be an adequate $\omega_1$-poset.
Let $\{\alpha_{\xi} : \xi < \omega_1 \}$ be a club subset of $\alpha$ with $\alpha_0 = 0$ and $\alpha_{\mu} < \alpha_{\xi}$ for $\mu < \xi < \omega_1$. For every countable successor ordinal $\xi = \mu + 1$, we define $\delta_{\xi} = \mbox{o.t.}(\al_{\xi}\setminus \al_{\mu})$. And we write $t_{\zeta} = \langle \omega_1 + 1,\zeta \rangle$ for $\zeta < \lambda$.  Let

$$T' = \un\{\{\xi\} \times \omega : \xi \mbox{ is a countable successor ordinal}\}\cup \{t_{\zeta} : \zeta \in \lambda\}.$$

\noindent
If $\xi = \mu + 1$  is a successor ordinal with $\mu \leq \omega_1$, for every $x\in \{\xi\}\times \omega$ we put $C_x = \{y\in \{\mu\}\times \omega : y \prec  x \}$. Then, for every countable successor ordinal $\xi = \mu + 1$ and every  $x\in \{\xi\}\times \omega$, we consider a compact Hausdorff scattered space $Z_x$ of height $\delta_{\xi} + 1$ such that $\mbox{CS}(Z_x) = \langle \omega \rangle_{\delta_{\xi}}$ and in such a way that $I_0(Z_x) = C_x$ and $I_{\delta_{\xi}}(Z_x) = \{x\}$. And for every ordinal $\zeta\in \lambda$ we consider a compact Hausdorff scattered space $Z_{t_{\zeta}}$ of height $\beta$ such that $\mbox{CS}(Z_{t_{\zeta}}) = \langle \omega \rangle_{\beta'}$ and in such a way that $I_0(Z_{t_{\zeta}}) = C_{t_{\zeta}}$ and  $I_{\be'}(Z_{t_{\zeta}}) = \{t_{\zeta}\}$. Also, we assume that for every $x\in T'$ we have  $Z_x \cap T = C_x \cup \{x\}$ and that for every $x,y \in T'$ with $x\neq y$ we have $(Z_x\setminus \{x\}) \cap (Z_y\setminus \{y\}) = \empt$.
Then, proceeding as in the proof of Theorem 2.5 (replacing in that proof  $S$ with $T$ and $S'$ with $T'$), we can construct an LCS space $Z$ whose underlying set is $T \cup \bigcup \{Z_x : x \in T'\}$ such that $\mbox{CS}(Z) = \langle \omega \rangle_{\alpha}\concat \langle
\lambda \rangle_{\beta}$.

\vspace{1mm}
If $\mbox{cf}(\al) = \om$, by using an adequate $\omega$-poset we can proceed by means of an argument similar to the one given in the preceding paragraph. And if $\al$ is a successor ordinal, we obtain the required LCS space by means of an adequate 1-poset. $\square$

\vspace{5mm}
Now, in order to complete the proof of Theorem 3.1, suppose that $V = L$ holds and $\lambda$ is an uncountable cardinal. Without loss of generality, we may assume that $\lambda\geq \omega_3$. We need to prove the following lemma.

\begin{lemma}
In some cardinal-preserving generic extension, there is an adequate $\om_1$-poset.
\end{lemma}

So, we will obtain the conclusion of Theorem 3.1 as an immediate consequence of Lemmas 3.6 and 3.7. Now, in order to show Lemma 3.7, note that  it follows from Theorem 3.3 that there is
a cardinal-preserving generic extension $N$ where there is an $\om_1$-strongly unbounded function $F : [\la]^2\longrightarrow \om_1$. Then, we construct an adequate $\om_1$-poset by means of a c.c.c. notion of forcing defined in $N$.

 We put $T = Y_{\om_1}$, $T_{\al} = \{\al\}\times \om$ for $\al < \om_1$, $T_{\om_1} = \{\om_1\}\times \lambda$ and $T_{\om_1 + 1} = \{\om_1 + 1\}\times \lambda$. As above, we write $t_{\xi} = \langle \omega_1 + 1,\xi \rangle$ for $\xi < \lambda$. Also, we put ${\Bbb B }_t = {\Bbb B }^{(\om_1)}_t$ for $t\in {\Bbb B}$.  Then, we define in $N$ the following poset $\pcal=\<P,\le\>$.  We say that  $p = \<X,\preceq,i\>\in P$ iff the following conditions hold:

\begin{enumerate}[(P1)]
\item $X\in [T]^{<\omega}$.
 \item $\preceq$ is a partial order on
$X$ such that
$x\prec y$ implies $\piz(x)<\piz (y)$.
 \item If $X\cap {\Bbb B }_{\zeta} \neq \emptyset$ for $\zeta \in\lambda$, then $t_{\zeta}\in X$.
\item If $x \in X\cap T_{\omega_1} \cap {\Bbb B }_{\zeta}$ for some $\zeta \in \lambda$, then $x \prec t_{\zeta}$ and $x \not\prec t_{\xi}$ for $\xi \neq \zeta$.
    \item If $s,t\in {\Bbb B}_t$ for some $t\in {\Bbb B}$ with $\pi(s) = \pi(t)$ and $s\neq t$, then $i\{s,t\} = \emptyset$.
\item If $t\in T_{\al + 1}$ for $\al \leq \om_1$ and $s\prec t$,  then there is a $v\in T_{\al}$ such that $s\preceq v \prec t$.
\item $\ifu:\br X;2; \to [X]^{<\omega}$ with $i\{x,y\} = \{x\}$ if $x \prec y$,  and such that if  $x,y \in X$  are $\preceq$-incomparable then the following conditions hold:

\begin{enumerate}[(a)]
\item \begin{equation}\notag
  \forall u\in X  ([u \preceq x\land u\preceq y]\text{ iff } u\preceq v \text{ for some } v \in i\{x,y\}).
  \end{equation}

\item If $x,y\in T_{\om_1}\cup T_{\om_1 + 1}$ with $\pi_B(x) \neq \pi_B(y)$, then $\pi[i\{x,y\}] \subset F\{\pi_B(x),\pi_B(y)\}$.

      \end{enumerate}
\end{enumerate}

The order on $P$ is the extension: $\<X',\preceq',\ifu'\>\le
\<X,\preceq,\ifu\>$
 iff $X \subset X'$, $\preceq=\preceq'\cap \, (X\times X)$ and
$\ifu\subs \ifu'$.

\begin{lemma}
(a) If $p = \langle X,\preceq,i\rangle\in P$ and $t\in T\setminus X$, then there is a $p' = \langle X',\preceq',i' \rangle\in P$ with $p'\leq p$ and $t\in X'$.

(b)  Assume that $p=\<X,\preceq,i\>\in  P$, $t\in X$,
$\alpha< \mbox{min}\{\pi(t),\om_1\}$ and $n <\omega$. Then,  there is a
$p'=\<X',\preceq',i'\>\in P$ with $p'\leq p$ and there is an $s\in
X'\setminus X$ with $\pi(s) = \al$ and $\rho(s)> n$ such that, for
every $x\in X$, $s\preceq' x$ iff $t\preceq x$.
\end{lemma}

\noindent {\bf Proof}. First, we prove $(a)$. Without loss of generality, we may assume that $t\in T_{\omega_1}$. Let $\zeta\in \lambda$ be such that $t\in {\Bbb B}_{\zeta}$. Assume that $X\cap {\Bbb B}_{\zeta} = \emptyset$. Otherwise, the argument is easier. We put $p' = \langle X',\preceq',i' \rangle$ where $X' = X \cup \{t,t_{\zeta}\}$, $\prec' = \prec \cup \{\langle s,t_{\zeta}\rangle : s\preceq v$ for some $v\in  T_{\omega_1}\cap {\Bbb  B}_{\zeta}\} \cup \{\langle t,t_{\zeta}\rangle\}$, $i'\{x,y\} = i\{x,y\}$ if $\{x,y\}\in [X]^2$, $i'\{x,t_{\zeta}\} = \{x\}$ if $x\preceq v$ for some $v\in T_{\omega_1}\cap {\Bbb  B}_{\zeta}$, $i'\{x,t_{\zeta}\} = \bigcup \{i\{x,v\} : v\in X\, \cap T_{\omega_1}\cap  {\Bbb  B}_{\zeta}\}$ if $x\in X$ and there is no $v\in T_{\omega_1}\cap {\Bbb  B}_{\zeta}$ such that $x\preceq v$, $i'\{x,t\} = \emptyset$ for every $x\in X$ and $i'\{t,t_{\zeta}\} = \{t\}$. It is easy to check that $p'$ is as required.

\vspace{1mm}
Now, we prove $(b)$. Let

$$L=\{{\alpha}\}\cup \{{\beta}: {\alpha} < {\be} < {\pi}(t) \land
\exists j<{\omega}\ {\beta}+j= \pi(t)\}.$$

\noindent Let ${\alpha}={\alpha}_0,\dots, {\alpha}_\ell$ be the increasing
enumeration of $L$.
Since $X$ is finite, for $j\leq l$ we can
pick an $s_j\in T_{\al_j}\setminus X$ such that  $\rho(s_j)> n$ if $\al_j < \om_1$ and $\pi_B(s_j) = \pi_B(t)$ if $\al_j = \om_1$ and $t\in T_{\om_1 + 1}$.
Let $X'=X\cup \{s_j:j\leq \ell\}$ and let
$$
\prec'=\prec \cup \{\langle s_j,y \rangle: t\preceq y\}\cup
\{\langle s_j,s_k \rangle :j<k\leq \ell\} .
$$
Now, we put $i'\{x,y\} = i\{x,y\}$ if $x,y\in [X]^2$, $i'\{s_j,s_k\}=\{s_j\}$ for $j < k \leq l$, $i'\{s_j,y\} =
\{s_j\}$ if $t\preceq y$
and $i'\{s_j,y\} = \emptyset$ otherwise.
It is easy to verify that  $\<X',\preceq',i'\>$ is as required. $\square$

\vspace{2mm} Assume that $\pcal=\<P,\le\>$ preserves cardinals. Let $G$ be a $\pcal$-generic filter. Put $ p = \langle X_p,\preceq_p,i_p\rangle$ for $p\in G$. By using Lemma 3.8, it is easy to see that $T = \un\{x_p : p\in G \}$ and if we put  $\preceq = \bigcup \{\preceq_p : p\in G \}$
and $i = \bigcup \{i_p : p\in G \}$, then  $\tcal = \<T,\preceq,i\>$ is an adequate $\om_1$-poset.  Then, in order to complete the proof of Lemma
3.7 we show the following lemma.

\begin{lemma}
$\pcal$ is c.c.c.
\end{lemma}

\noindent {\bf Proof}. Assume that $ R = \<r_{\nu}:{\nu}<\omega_1\>\subs P$ with $r_{\nu}\neq r_{\mu}$ for $\nu <\mu<\omega_1$. For $\nu < \omega_1$, we write
 $r_{\nu}=\<X_{\nu},\preceq_{\nu},\ifu_{\nu}\>$. By thinning
out $\<r_{\nu}:{\nu}<\omega_1\>$ by  standard combinatorial
arguments, we can assume the following:

\begin{enumerate}[(A)]
 \item \begin{enumerate}[(a)]
\item $\{X_{\nu} : \nu < \omega_1 \}$ forms a $\Delta$-system with kernel $X$,
\item $|X_{\nu}| = |X_{\mu}|$ for $\nu < \mu < \om_1$,
\item $\pi[X\cap \blocks_S]$ is an initial segment of $\pi[X_{\nu}]$ for each $\nu < \omega_1$,
\item for each $\alpha < \omega_1$, either $X_{\nu}\cap (\{\alpha\} \times \omega) = X \cap (\{\alpha\} \times \omega)$ for every $\nu < \omega_1$ or there is at most one $\nu < \omega_1$ such that $X_{\nu}\cap (\{\alpha\} \times \omega) \neq \empt$.
\end{enumerate}

\item If $H_{\nu} = \{t\in {\Bbb B} : X_{\nu}\cap {\Bbb B}_t \neq \emptyset \}$ for each $\nu < \om_1$, then $\{H_{\nu} : \nu < \om_1 \}$ forms a $\Delta$-system such that the following conditions hold:

    \begin{enumerate}[(a)]
    \item $|H_{\nu}| = |H_{\mu}|$ for $\nu < \mu < \om_1$,

    \item for every $\zeta \in \lambda$, if $X\cap {\Bbb B}_{\zeta} \neq \emptyset$ then $X_{\nu} \cap {\Bbb B}_{\zeta} = X \cap {\Bbb B}_{\zeta}$ for each $\nu < \om_1$.

        \end{enumerate}

\item For each ${\nu}<{\mu}<\omega_1$ there is an isomorphism
$h=h_{{\nu},{\mu}}:\<X_{\nu},\preceq_{\nu}, \ifu_{\nu}\>\to
\<X_{\mu},\preceq_{\mu},\ifu_{\mu}\>$ such that:
\begin{enumerate}[(a)]
  \item $h\restriction X=\operatorname{id},$
  \item $\pib(x)=S$ iff $\pib(h(x))=S,$
\item $\pi(x) = \om_1$ iff $\pi(h(x)) = \om_1$,
\item $\pi(x) = \om_1 + 1$ iff $\pi(h(x)) = \om_1 + 1$,
\item $\pib(x)=\pib(y)$  iff $\pib(h(x))=\pib(h(y))$,
\item $\ifu_{\nu}\{x,y\}=\ifu_{\mu}\{x,y\}$ for all $\{x,y\}\in [X]^2$.

\end{enumerate}

\end{enumerate}

\vspace{1mm} To verify condition $(C)(f)$, assume that $x,y\in X$ are $\preceq_{\nu}$-incomparable for $\nu < \omega_1$. If either $x\in {\Bbb B}_S$ or $y\in {\Bbb B}_S$, the case is obvious. So, assume that $x,y\in T_{\om_1}\cup T_{\om_1 + 1}$. If $\pi_B(x) = \pi_B(y)$, by conditions $(P4)$ and $(P5)$, we deduce that $i_{\nu}(x,y\} = \emptyset$ for $\nu < \om_1$. And if $\pi_B(x) \neq \pi_B(y)$, we apply condition $(P7)(b)
$ to obtain $i_{\nu}\{x,y\} = i_{\mu}\{x,y\}$ for $\nu < \mu < \om_1$.

\vspace{1mm}
Let $\delta = \mbox{max}(\pi[X\cap {\Bbb B}_S])$. Since $F : [\la]^2\longrightarrow \om_1$ is an $\om_1$-strongly unbounded function, we deduce from condition $(B)$ that there are ordinals $\nu <\mu < \om_1$ such that if we put $D = \{\xi\in\lambda :  (X_{\nu}\setminus X) \cap {\Bbb B}_{\xi} \neq \emptyset\}$ and $E = \{\xi\in\lambda :  (X_{\mu}\setminus X) \cap {\Bbb B}_{\xi} \neq \emptyset\}  $, then $F\{\xi,\zeta\} > \delta$ for every $\xi\in D$ and every $\zeta\in E$. We show that $r_{\nu}$ and $r_{\mu}$ are
compatible in $\pcal$.  We put $p= r_{\nu}$ and $q = r_{\mu}$. And we write $p=\<X_p,\preceq_p,i_p\>$ and $q=\<X_q,\preceq_q,i_q\>$.
Then, we define the extension $r = \<X_r,\preceq_r,i_r\>$ of $p$
and $q$ as follows. We put $X_r = X_p\cup X_q$ and $\leq_r = \leq_p \cup \leq_q$. Clearly, $\leq_r$ is a partial order on $X_r$. So, we define the infimum function $i_r$. We put $i_r\{s,t\} = i_p\{s,t\}$ if $\{s,t\}\in [X_p]^2$, and $i_r\{s,t\} = i_q\{s,t\}$ if $\{s,t\}\in [X_q]^2$. By condition $(C)(f)$, $i_r \upharpoonright [X]^2$ is well-defined. Now, assume that $s,t\in X_r$ with $s\in X_p\setminus X_q$ and $t\in X_q\setminus X_p$. Note that $s,t$ are not comparable in $\langle X_r,\leq_r\rangle$ and there is no $u\in (X_p\cup X_q)\setminus X$ such that $u \preceq_r s,t$. Then, we define $i_r\{s,t\} = \{ u\in X \cap {\Bbb B}_S : u \prec_r s,t \}$. It is easy to check that $r\in P$, and so $r\leq p,q$.  $\square$

\vspace{1mm} This concludes the proofs of Lemma 3.9, Lemma 3.7 and Theorem 3.1.

 \vspace{6mm}{\bf Acknowledgements.} The first author was supported by the Spanish Ministry of Education DGI grant MTM2017-86777-P and by the Catalan DURSI grant 2017SGR270. The second author was supported by NKFIH grant nos. K113047 and K129211.

\vspace{6mm}

FACULTAT  DE MATEM\`ATIQUES I INFORM\`ATICA, UNIVERSITAT DE BARCELONA, GRAN VIA 585, 08007 BARCELONA, SPAIN

\vspace{1mm} {\em E-mail address}: jcmartinez@ub.edu

\vspace{3mm}

ALFR\'ED R\'ENYI INSTITUTE OF MATHEMATICS

\vspace{1mm} {\em E-mail address}: soukup@renyi.hu


\begin{thebibliography}{9}


   \bibitem{ba} J. Bagaria, Thin-tall spaces and cardinal sequences, {\bf Open problems in Topology II} (E. Peral, editor), Elsevier, Amsterdam, 2007, pp. 115-124.


\bibitem{bs} J. E. Baumgartner and S. Shelah, {\em Remarks on
superatomic Boolean algebras}, {\bf Annals of Pure and  Applied Logic}, vol. 33 (l987), no. 2, pp. 109-129.



\bibitem{rv} K. Er-Rhaimini and B. Veli{\v c}kovi{\' c}, {\em PCF structures of height less than $\om_3$}, {\bf The Journal of Symbolic Logic}, vol 75 (2010), no. 4, pp. 1231-1248.


\bibitem{jsw} I. Juh\'asz, L. Soukup and W. Weiss, {\em  Cardinal sequences of length $< \omega_2$ under GCH}, {\bf Fundamenta Mathematicae}, vol. 189 (2006), no. 1, pp. 35-52.

\bibitem{jw} I. Juh\'asz and W. Weiss, {\em Cardinal sequences},  {\bf Annals of Pure and Applied Logic}, vol.  144 (2006), no. 1-3, pp. 96-106.

\bibitem{ko} P. Koszmider, {\em  Universal matrices and strongly unbounded functions}, {\bf Mathematical Research Letters}, vol. 9 (2002), no. 4,  pp. 549-566.

\bibitem{ma1} J. C. Mart\'{\i}nez, {\em A consistency result on thin-very tall Boolean algebras}, {\bf Israel Journal of Mathematics}, vol. 123 (2001), pp. 273-284.

\bibitem{ma2} J. C. Mart\'{\i}nez, Cardinal sequences for superatomic Boolean algebras, {\bf Infinity, computability and metamathematics} (S. Geschke, B. L\"{o}we and P. Schlicht, editors) Tributes Series, vol. 23, College Publications, Milton Keynes, 2014,  pp. 273-284.

 \bibitem{ms} J. C. Mart\'{\i}nez and L. Soukup, {\em Superatomic
Boolean algebras constructed from strongly unbounded functions}, {\bf Mathematical Logic Quarterly}, vol. 57 (2011), no. 5, pp. 456-469.

\bibitem{ro} J. Roitman, {\em A very thin thick superatomic Boolean algebra}, {\bf Algebra Universalis}, vol.  21 (1985), no. 2-3, pp. 137-142.

\bibitem{so} L. Soukup, {\em Wide scattered spaces and morasses}, {\bf Topology and its Applications}, vol. 158 (2011), no. 5, pp. 697-707.

\end{thebibliography}
\end{document}